\pdfoutput=1
\RequirePackage{silence}
\documentclass[a4paper,svgnames]{amsart}
\usepackage[hmarginratio={1:1},vmarginratio={1:1},lmargin=60.0pt,tmargin=60.0pt]{geometry}

\usepackage{booktabs}
\allowdisplaybreaks

\usepackage[T1]{fontenc}
\usepackage{latexsym,exscale,amsfonts,amssymb,mathtools}
\usepackage{amsmath,amsthm,amsfonts,amssymb,amscd,textcomp,bbm}
\usepackage{mathrsfs,stackrel}
\usepackage{etoolbox}
\usepackage{tikz,tikz-cd}
\usetikzlibrary{matrix,arrows.meta,decorations.markings,shapes.multipart,decorations.pathreplacing,decorations.pathmorphing,shapes.geometric}
\usepackage{aliascnt}
\usepackage{ytableau}
\usepackage{xparse}
\usepackage{cite}

\usepackage{dynkin-diagrams}
\usepackage{nicematrix}
\usepackage{tikz-3dplot}
\usepackage{mathrsfs}
\DeclareMathAlphabet{\mathscrbf}{OMS}{mdugm}{b}{n}

\usepackage{array}
\newcolumntype{C}{>{$}c<{$}}

\usepackage{xcolor,colortbl}

\definecolor{mygray}{gray}{0.6}
\definecolor{mygraydark}{gray}{0.4}
\definecolor{mygraylight}{gray}{0.85}
\definecolor{spinach}{RGB}{46,139,87}
\definecolor{tomato}{RGB}{255,99,71}
\definecolor{orchid}{RGB}{143,40,194}
\definecolor{neon}{RGB}{77,77,255}
\definecolor{lightneon}{RGB}{110,110,255}
\definecolor{pumpkin}{RGB}{224,180,80}
\definecolor{citron}{RGB}{190,180,90}

\definecolor{lava}{RGB}{207,16,32}
\definecolor{cream}{RGB}{255,253,208}
\definecolor{verdigris}{RGB}{67,179,174}
\definecolor{Black}{RGB}{0,0,0}
\definecolor{mydarkblue}{RGB}{10,10,170}
\definecolor{darkspinach}{RGB}{20,70,20}
\definecolor{darktomato}{RGB}{155,40,30}
\definecolor{darkorchid}{RGB}{50,10,100}
\definecolor{darklava}{RGB}{150,8,16}

\definecolor{zero}{RGB}{0,0,0}
\definecolor{one}{RGB}{255,0,0}
\definecolor{two}{RGB}{0,255,0}
\definecolor{three}{RGB}{0,0,255}

\usepackage{todonotes}

\def\changed#1{{#1}}
\def\ochanged#1{{#1}}

\usepackage{enumitem}
\setlist[enumerate]{itemsep=0.15cm,label=\emph{\upshape(\alph*)}}
\setlist[enumerate,2]{itemsep=0.15cm,label=\emph{\upshape(\roman*)}}
\setlist[enumerate,3]{itemsep=0.15cm,label=\emph{\upshape(\Alph*)}}

\let\emph\relax
\DeclareTextFontCommand{\emph}{\bfseries\em}

\renewcommand{\dots}{\text{...}}

\newcommand{\ie}{\text{i.e.}}
\newcommand{\eg}{\text{e.g.}}

\newcommand{\etc}{\text{etc.}}

\usepackage[all]{xy}
\usepackage{tikz}
\usetikzlibrary{cd}
\usetikzlibrary{decorations}
\usetikzlibrary{decorations.markings}
\usetikzlibrary{decorations.pathreplacing}
\usetikzlibrary{decorations.pathmorphing}
\usetikzlibrary{arrows.meta,shapes,positioning,matrix,calc}
\usetikzlibrary{shapes.callouts}
\usetikzlibrary{tqft}
\tikzset{anchorbase/.style={baseline={([yshift=-0.5ex]current bounding box.center)}},
tinynodes/.style={font=\tiny,text height=0.25ex,text depth=0.05ex},
smallnodes/.style={font=\scriptsize,text height=0.75ex,text depth=0.15ex},
usual/.style={line width=2.0,color=black},
dot/.style = {
decoration={markings,
post length=0.25mm,
pre length=0.25mm,
mark=at position #1 with {\node[circle,radius=0.2cm,inner sep=-1.5pt,color=black,fill=black]{};}
},
postaction={decorate}
},
dot/.default=1,
dotg/.style = {
decoration={markings,
post length=0.25mm,
pre length=0.25mm,
mark=at position #1 with {\node[circle,radius=0.5cm,inner sep=-2.5pt,color=neon,fill=neon]{};}
},
postaction={decorate}
},
dotg/.default=0.5,
dotw/.style = {
decoration={markings,
post length=0.25mm,
pre length=0.25mm,
mark=at position #1 with {\node[circle,radius=0.2cm,inner sep=-1.5pt,color=black,fill=citron]{};}
},
postaction={decorate}
},
dotw/.default=1,
}
\tikzstyle directed=[postaction={decorate,decoration={markings,mark=at position #1 with {\arrow[line width=0.5mm, black]{>}}}}]
\tikzstyle rdirected=[postaction={decorate,decoration={markings,mark=at position #1 with {\arrow[line width=0.5mm, black]{<}}}}]

\newcommand{\C}{\mathbb{C}}
\newcommand{\R}{\mathbb{R}}

\newcommand{\N}{\mathbb{Z}_{\geq 0}}

\newcommand{\Z}{\mathbb{Z}}

\newcommand{\K}{\mathbb{K}}
\newcommand{\F}{\mathbb{F}}

\newcommand{\setstuff}[1]{\mathrm{#1}}
\newcommand{\catstuff}[1]{\mathbf{#1}}

\newcommand{\obstuff}[1]{\mathtt{#1}}

\newcommand{\End}{\setstuff{End}}
\newcommand{\Hom}{\setstuff{Hom}}

\newcommand{\cob}[1][k]{\mathcal{C}\mathrm{ob}_{#1}}
\newcommand{\sym}[1][n]{\mathcal{S}_{#1}}
\newcommand{\tlmon}[1][n]{\mathcal{TL}_{#1}}
\newcommand{\brmon}[1][n]{\mathcal{B}\mathrm{r}_{#1}}
\newcommand{\pamon}[1][n]{\mathcal{P}\mathrm{a}_{#1}}
\newcommand{\robrmon}[1][n]{\mathcal{R}\mathrm{o}\mathcal{B}\mathrm{r}_{#1}}
\newcommand{\romon}[1][n]{\mathcal{R}\mathrm{o}_{#1}}
\newcommand{\ppamon}[1][n]{\mathrm{p}\mathcal{P}\mathrm{a}_{#1}}
\newcommand{\probrmon}[1][n]{\mathrm{p}\mathcal{R}\mathrm{o}\mathcal{B}\mathrm{r}_{#1}}
\newcommand{\momon}[1][n]{\mathcal{M}\mathrm{o}_{#1}}

\newcommand{\promon}[1][n]{\mathrm{p}\mathcal{R}\mathrm{o}_{#1}}
\newcommand{\pbrmon}[1][n]{\mathrm{p}\mathcal{B}\mathrm{r}_{#1}}
\newcommand{\psym}[1][n]{\mathrm{p}\mathcal{S}_{#1}}
\newcommand{\obrmon}[1][n]{\mathrm{o}\mathcal{B}\mathrm{r}_{#1}}
\newcommand{\Rep}{\mathrm{Rep}}
\newcommand{\GL}{\mathrm{GL}}

\newcommand{\repsnn}[1][n]{\Rep(\mathrm{S}_{#1})}
\newcommand{\repgll}[1][n]{\Rep\big(\mathrm{GL}_{#1}(\C)\big)}
\newcommand{\reposp}[1][n]{\Rep\big(\mathrm{OSP}_{#1}(\C)\big)}
\newcommand{\repospt}[1][n]{\Rep\big(\mathrm{O}_{#1}(\C)\big)}
\newcommand{\repgl}[1][n]{\Rep\big(\mathrm{GL}_{#1}(\mathbb{F}_{q})\big)}

\def\NewTheorem#1{%
\newaliascnt{#1}{equation}%
\newtheorem{#1}[#1]{#1}%
\aliascntresetthe{#1}%
\expandafter\def\csname #1autorefname\endcsname{#1}%
}
\def\equationautorefname~#1\null{(#1)\null}

\numberwithin{equation}{subsection}

\NewTheorem{Proposition}
\NewTheorem{Theorem}
\NewTheorem{Corollary}
\AtEndEnvironment{Corollary}{\null\hfill$\square$}%
\NewTheorem{Lemma}
\NewTheorem{Conjecture}
\NewTheorem{Speculation}
\NewTheorem{Question}
\NewTheorem{Assumption}
\theoremstyle{definition}
\NewTheorem{Definition}
\AtEndEnvironment{Definition}{\null\hfill$\Diamond$}%
\NewTheorem{Classification Problem}
\AtEndEnvironment{Classification Problem}{\null\hfill$\Diamond$}%
\NewTheorem{Notation}
\AtEndEnvironment{Notation}{\null\hfill$\Diamond$}%
\NewTheorem{Example}
\AtEndEnvironment{Example}{\null\hfill$\Diamond$}%
\NewTheorem{Examples}
\AtEndEnvironment{Examples}{\vskip-10mm\null\hfill$\Diamond$}%

\theoremstyle{remark}
\NewTheorem{Remark}
\AtEndEnvironment{Remark}{\null\hfill$\Diamond$}%

\usepackage[hypertexnames=false]{hyperref}
\usepackage{bookmark}
\hypersetup{
pdftoolbar=true,
pdfmenubar=true,
pdffitwindow=false,
pdfstartview={FitH},
pdftitle={Growth problems in diagram categories},
pdfauthor={Jonathan Gruber and Daniel Tubbenhauer},
pdfsubject={},
pdfcreator={Jonathan Gruber and Daniel Tubbenhauer},
pdfproducer={Jonathan Gruber and Daniel Tubbenhauer},
pdfkeywords={},
pdfnewwindow=true,
colorlinks=true,
linkcolor=mydarkblue,
citecolor=teal,
filecolor=magenta,
urlcolor=orchid,
linkbordercolor=lava,
citebordercolor=teal,
urlbordercolor=orchid,
linktocpage=true
}

\def\makeautorefname#1#2{\csdef{#1autorefname}{#2}}
\makeautorefname{section}{Section}%
\makeautorefname{subsection}{Section}%
\makeautorefname{subsubsection}{Section}%

\begin{document}
\title[Growth problems in diagram categories]{Growth problems in diagram categories}
\author[J. Gruber and D. Tubbenhauer]{Jonathan Gruber and Daniel Tubbenhauer}

\address{J.G.: FAU Erlangen-Nuremberg, Department of Mathematics, Cauerstr.\ 11, 91058 Erlangen, Germany, \href{https://jgruber4.github.io/}{jgruber4.github.io/}, \href{https://orcid.org/0000-0001-5975-8041}{ORCID 0000-0001-5975-8041}}
\email{jonathan.gruber@fau.de}

\address{D.T.: The University of Sydney, School of Mathematics and Statistics F07, Office Carslaw 827, NSW 2006, Australia, \href{http://www.dtubbenhauer.com}{www.dtubbenhauer.com}, \href{https://orcid.org/0000-0001-7265-5047}{ORCID 0000-0001-7265-5047}}
\email{daniel.tubbenhauer@sydney.edu.au}

\begin{abstract}
In the semisimple case, we derive (asymptotic) formulas for the growth rate of the number of summands in tensor powers of the generating object in diagram/interpolation categories.
\end{abstract}


\addtocontents{toc}{\protect\setcounter{tocdepth}{1}}

\maketitle

\tableofcontents

\arrayrulewidth=0.5mm
\setlength{\arrayrulewidth}{0.5mm}

\section{Introduction}\label{S:Intro}

We begin with the following table, the meaning of which we will explain shortly:
\begin{gather}\label{Eq:Main}
\begin{gathered}
\begin{tabular}{c|c|c|c}
\arrayrulecolor{tomato}
Symbol & Diagrams & $\obstuff{X}$ & $\sqrt[n]{b_{n}}$ $\sim$
\\
\hline
\hline
$\cob$ & \begin{tikzpicture}[anchorbase]
\draw[usual] (0.5,0) to[out=90,in=180] (1.25,0.45) to[out=0,in=90] (2,0);
\draw[usual] (0.5,0) to[out=90,in=180] (1,0.35) to[out=0,in=90] (1.5,0);
\draw[usual] (0,1) to[out=270,in=180] (0.75,0.55) to[out=0,in=270] (1.5,1);
\draw[usual] (1.5,1) to[out=270,in=180] (2,0.55) to[out=0,in=270] (2.5,1);
\draw[usual,dotg=0.25,dotg=0.85] (0,0) to (0.5,1);
\draw[usual] (1,0) to (1,1);
\draw[usual,dotg] (2.5,0) to (2.5,1);
\draw[usual,dot] (2,1) to (2,0.8);
\end{tikzpicture} & $\bullet$ & $\tfrac{2}{e}\cdot\tfrac{n}{\log 2n/k}$
\end{tabular}
\\
\begin{tabular}{c|c|c|c||c|c|c|c}
\arrayrulecolor{tomato}
Symbol & Diagrams & $\obstuff{X}$ & $\sqrt[n]{b_{n}}$ $\sim$
& Symbol & Diagrams & $\obstuff{X}$ & $\sqrt[n]{b_{n}}$ $\sim$
\\
\hline
\hline
$\ppamon[t]$ & \begin{tikzpicture}[anchorbase]
\draw[usual] (0.5,0) to[out=90,in=180] (1.25,0.45) to[out=0,in=90] (2,0);
\draw[usual] (0.5,0) to[out=90,in=180] (1,0.35) to[out=0,in=90] (1.5,0);
\draw[usual] (0.5,1) to[out=270,in=180] (1,0.55) to[out=0,in=270] (1.5,1);
\draw[usual] (1.5,1) to[out=270,in=180] (2,0.55) to[out=0,in=270] (2.5,1);
\draw[usual] (0,0) to (0,1);
\draw[usual] (2.5,0) to (2.5,1);
\draw[usual,dot] (1,0) to (1,0.2);
\draw[usual,dot] (1,1) to (1,0.8);
\draw[usual,dot] (2,1) to (2,0.8);
\end{tikzpicture} & $\bullet$ & $4$
& $\pamon[t]$ & \begin{tikzpicture}[anchorbase]
\draw[usual] (0.5,0) to[out=90,in=180] (1.25,0.45) to[out=0,in=90] (2,0);
\draw[usual] (0.5,0) to[out=90,in=180] (1,0.35) to[out=0,in=90] (1.5,0);
\draw[usual] (0,1) to[out=270,in=180] (0.75,0.55) to[out=0,in=270] (1.5,1);
\draw[usual] (1.5,1) to[out=270,in=180] (2,0.55) to[out=0,in=270] (2.5,1);
\draw[usual] (0,0) to (0.5,1);
\draw[usual] (1,0) to (1,1);
\draw[usual] (2.5,0) to (2.5,1);
\draw[usual,dot] (2,1) to (2,0.8);
\end{tikzpicture} & $\bullet$ & $\tfrac{2}{e}\cdot\tfrac{n}{\log 2n}$
\\
\hline
$\momon[t]$ & \begin{tikzpicture}[anchorbase]
\draw[usual] (0.5,0) to[out=90,in=180] (1.25,0.5) to[out=0,in=90] (2,0);
\draw[usual] (1,0) to[out=90,in=180] (1.25,0.25) to[out=0,in=90] (1.5,0);
\draw[usual] (2,1) to[out=270,in=180] (2.25,0.75) to[out=0,in=270] (2.5,1);
\draw[usual] (0,0) to (1,1);
\draw[usual,dot] (2.5,0) to (2.5,0.2);
\draw[usual,dot] (0,1) to (0,0.8);
\draw[usual,dot] (0.5,1) to (0.5,0.8);
\draw[usual,dot] (1.5,1) to (1.5,0.8);
\end{tikzpicture} & $\bullet$ & $3$
& $\robrmon[t]$ & \begin{tikzpicture}[anchorbase]
\draw[usual] (1,0) to[out=90,in=180] (1.25,0.25) to[out=0,in=90] (1.5,0);
\draw[usual] (1,1) to[out=270,in=180] (1.75,0.55) to[out=0,in=270] (2.5,1);
\draw[usual] (0,0) to (0.5,1);
\draw[usual] (2.5,0) to (2,1);
\draw[usual,dot] (0.5,0) to (0.5,0.2);
\draw[usual,dot] (2,0) to (2,0.2);
\draw[usual,dot] (0,1) to (0,0.8);
\draw[usual,dot] (1.5,1) to (1.5,0.8);
\end{tikzpicture} & $\bullet$ & $\tfrac{\sqrt{2}}{\sqrt{e}}\cdot\sqrt{n}$
\\
\hline
$\tlmon[t]$ & \begin{tikzpicture}[anchorbase]
\draw[usual] (0.5,0) to[out=90,in=180] (1.25,0.5) to[out=0,in=90] (2,0);
\draw[usual] (1,0) to[out=90,in=180] (1.25,0.25) to[out=0,in=90] (1.5,0);
\draw[usual] (0,1) to[out=270,in=180] (0.25,0.75) to[out=0,in=270] (0.5,1);
\draw[usual] (2,1) to[out=270,in=180] (2.25,0.75) to[out=0,in=270] (2.5,1);
\draw[usual] (0,0) to (1,1);
\draw[usual] (2.5,0) to (1.5,1);
\end{tikzpicture} & $\bullet$ & $2$
& $\brmon[t]$ & \begin{tikzpicture}[anchorbase]
\draw[usual] (0.5,0) to[out=90,in=180] (1.25,0.45) to[out=0,in=90] (2,0);
\draw[usual] (1,0) to[out=90,in=180] (1.25,0.25) to[out=0,in=90] (1.5,0);
\draw[usual] (0,1) to[out=270,in=180] (0.75,0.55) to[out=0,in=270] (1.5,1);
\draw[usual] (1,1) to[out=270,in=180] (1.75,0.55) to[out=0,in=270] (2.5,1);
\draw[usual] (0,0) to (0.5,1);
\draw[usual] (2.5,0) to (2,1);
\end{tikzpicture} & $\bullet$ & $\tfrac{\sqrt{2}}{\sqrt{e}}\cdot\sqrt{n}$
\\
\hline
$\promon[t]$ & \begin{tikzpicture}[anchorbase]
\draw[usual] (0,0) to (0.5,1);
\draw[usual] (0.5,0) to (1,1);
\draw[usual] (2,0) to (1.5,1);
\draw[usual] (2.5,0) to (2.5,1);
\draw[usual,dot] (1,0) to (1,0.2);
\draw[usual,dot] (1.5,0) to (1.5,0.2);
\draw[usual,dot] (0,1) to (0,0.8);
\draw[usual,dot] (2,1) to (2,0.8);
\end{tikzpicture} & $\bullet$ & $2$
& $\romon[t]$ & \begin{tikzpicture}[anchorbase]
\draw[usual] (0,0) to (1,1);
\draw[usual] (0.5,0) to (0,1);
\draw[usual] (2,0) to (2,1);
\draw[usual] (2.5,0) to (0.5,1);
\draw[usual,dot] (1,0) to (1,0.2);
\draw[usual,dot] (1.5,0) to (1.5,0.2);
\draw[usual,dot] (1.5,1) to (1.5,0.8);
\draw[usual,dot] (2.5,1) to (2.5,0.8);
\end{tikzpicture} & $\bullet$ & $\tfrac{\sqrt{1}}{\sqrt{e}}\cdot\sqrt{n}$
\\
\hline
$\psym[t]$ & \begin{tikzpicture}[anchorbase]
\draw[usual] (0,0) to (0,1);
\draw[usual] (0.5,0) to (0.5,1);
\draw[usual] (1,0) to (1,1);
\draw[usual] (1.5,0) to (1.5,1);
\draw[usual] (2,0) to (2,1);
\draw[usual] (2.5,0) to (2.5,1);
\end{tikzpicture} & $\bullet$ & $1$
& $\sym[t]$ & \begin{tikzpicture}[anchorbase]
\draw[usual] (0,0) to (1,1);
\draw[usual] (0.5,0) to (0,1);
\draw[usual] (1,0) to (1.5,1);
\draw[usual] (1.5,0) to (2.5,1);
\draw[usual] (2,0) to (2,1);
\draw[usual] (2.5,0) to (0.5,1);
\end{tikzpicture} & $\bullet$ & $\tfrac{\sqrt{1}}{\sqrt{e}}\cdot\sqrt{n}$
\end{tabular}
\\
\begin{tabular}{c|c|c|c||c|c|c|c}
\arrayrulecolor{tomato}
Symbol & Diagrams & $\obstuff{X}$ & $\sqrt[n]{b_{n}}$ $\sim$
& Symbol & Diagrams & $\obstuff{X}$ & $\sqrt[n]{b_{n}}$ $\sim$
\\
\hline
\hline
$\obrmon[t]$ & \begin{tikzpicture}[anchorbase]
\draw[usual,directed=1] (0,0) to (1,1);
\draw[usual,directed=1] (0.5,0) to (0,1);
\draw[usual,directed=1] (1,0) to (1.5,1);
\draw[usual,directed=1] (1.5,0) to (2.5,1);
\draw[usual,directed=1] (2,0) to (2,1);
\draw[usual,directed=1] (2.5,0) to (0.5,1);
\end{tikzpicture} & $\uparrow$ & $\tfrac{\sqrt{1}}{\sqrt{e}}\cdot\sqrt{n}$
& 
$\obrmon[t]$ & \begin{tikzpicture}[anchorbase]
\draw[usual,directed=1] (0.5,0) to[out=90,in=180] (1.25,0.45) to[out=0,in=90] (2,0);
\draw[usual,directed=1] (1.5,0) to[out=90,in=0] (1.25,0.25) to[out=180,in=90] (1,0);
\draw[usual,directed=1] (0,1) to[out=270,in=180] (0.75,0.55) to[out=0,in=270] (1.5,1);
\draw[usual,directed=1] (1,1) to[out=270,in=180] (1.75,0.55) to[out=0,in=270] (2.5,1);
\draw[usual,directed=1] (0,0) to (0.5,1);
\draw[usual,directed=1] (2,1) to (2.5,0);
\end{tikzpicture} & $\uparrow\downarrow$ & $\tfrac{2}{e}\cdot n$
\\
\end{tabular}
\\
\begin{tabular}{c|c|c|c}
\arrayrulecolor{tomato}
Symbol & Diagrams & $\obstuff{X}$ & $\sqrt[n]{b_{n}}$ $\sim$
\\
\hline
\hline
$\repgl[t]$ & \begin{tikzpicture}[anchorbase]
\draw[usual] (0.5,0) to[out=90,in=180] (1.25,0.45) to[out=0,in=90] (2,0);
\draw[usual] (0.5,0) to[out=90,in=180] (1,0.35) to[out=0,in=90] (1.5,0);
\draw[usual] (0,1) to[out=270,in=180] (0.75,0.55) to[out=0,in=270] (1.5,1);
\draw[usual] (1.5,1) to[out=270,in=180] (2,0.55) to[out=0,in=270] (2.5,1);
\draw[usual] (0,0) to (0.5,1);
\draw[usual] (1,0) to (1,1);
\draw[usual] (2.5,0) to (2.5,1);
\draw[usual,dot] (2,1) to (2,0.8);
\end{tikzpicture}
+
\begin{tikzpicture}[anchorbase]
\draw[white] (0,0) to (0,1);
\draw[usual,dotw] (0,0) to (0,0.5);
\end{tikzpicture}
+
\begin{tikzpicture}[anchorbase]
\draw[usual] (0,0) to (0.5,0.5);
\draw[usual] (1,0) to (0.5,0.5);
\draw[usual] (0.5,0.5) to (0.5,1);
\node[circle, draw,fill=white,inner sep=0.0cm] at (0.5,0.5) {$+$};
\end{tikzpicture}
+
\begin{tikzpicture}[anchorbase]
\draw[usual] (0,0) to (0,1);
\node[circle, draw,fill=white,inner sep=0.02cm] at (0,0.5) {$a$};
\end{tikzpicture}
for $a\in\mathbb{F}_{q}$ & $\bullet$ & $q^{(n+1)/2}$
\\
\end{tabular}
\end{gathered}
\end{gather}
Let $\catstuff{C}$ be an
additive Krull--Schmidt monoidal category.
Let $\obstuff{X}\in\catstuff{C}$ be an object of $\catstuff{C}$.
We define 
\begin{gather*}
b_{n}=b_{n}^{\catstuff{C},\obstuff{X}}:=\#\text{indecomposable summands in $\obstuff{X}^{\otimes n}$ counted with multiplicities}.
\end{gather*}
When $\catstuff{C}$ is semisimple, then $b_{n}=l_{n}$, the latter counting the number of simple factors in $\obstuff{X}^{\otimes n}$, and we will use this implicitly throughout.

\begin{Notation}
We also consider the function $\N\to\N,n\mapsto b_{n}$, and the corresponding sequence $(b_n)_{n\in\N}$, both denoted by the same symbol $b_n$. More generally, we identify sequences with their associated functions and values.
\end{Notation}

The function $b_{n}$ has been the subject of extensive study; see, for example, \cite{CoEtOsTu-growth-fractal,LaPoReSo-growth-qgroup,La-char2-story,He-1,He-2} for some recent work. In particular, in well-behaved categories, such as \ochanged{finite dimensional} representations of a group \cite{CoOsTu-growth,LaTuVa-growth-pfdim,LaTuVa-growth-pfdim-inf}, one has
\begin{gather*}
\lim_{n\to\infty}\sqrt[n]{b_{n}}\in\R_{\geq 1},
\quad\text{exponential growth}
\end{gather*}
which shows that $b_{n}$ grows \emph{exponentially}. In contrast, some still well-structured categories exhibit \emph{superexponential} growth, meaning that 
\begin{gather*}
\text{$\sqrt[n]{b_{n}}$ is unbounded,}
\quad\text{superexponential growth}
\end{gather*}
as observed\ochanged{, for example,} in \cite{De-cat-st}.

We study the asymptotic behavior of $\sqrt[n]{b_{n}}$ in the following cases, all in the semisimple situation (all parameters are generic) and over the complex numbers. Our starting point are prototypical examples of \emph{diagram categories}: quotients of the \emph{cobordism category} $\cob$, see \cite{KhSa-cobordisms,KhOsKo-cobordisms} \ochanged{(or \cite{Ko-tqfts} for the purely topological incarnation)}. Here the objects are one-dimensional compact manifolds (circles) and the morphisms are two-dimensional cobordisms (pants), which we will draw using their spines with handles as dots:
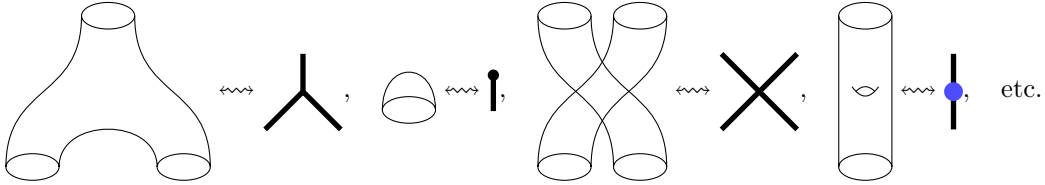
\begin{gather*}
\begin{tikzpicture}[tqft/cobordism/.style={draw},tqft/every lower boundary component/.style={draw=black},anchorbase]
\pic[tqft/pair of pants];
\end{tikzpicture}
\leftrightsquigarrow
\begin{tikzpicture}[anchorbase]
\draw[usual] (0,0) to (0.5,0.5) to (1,0);
\draw[usual] (0.5,0.5) to (0.5,1);
\end{tikzpicture}
,\quad
\begin{tikzpicture}[tqft/cobordism/.style={draw},tqft/every lower boundary component/.style={draw=black},anchorbase]
\pic[tqft/cap];
\end{tikzpicture}
\leftrightsquigarrow
\begin{tikzpicture}[anchorbase]
\draw[usual,dot] (0,0) to (0,0.5);
\end{tikzpicture}
,\quad
\begin{tikzpicture}[tqft/cobordism/.style={draw},tqft/every lower boundary component/.style={draw=black},anchorbase]
\pic[tqft/cylinder to prior, at={(0,0)}];
\pic[tqft/cylinder to next, at={(-1,0)}];
\end{tikzpicture}
\leftrightsquigarrow
\begin{tikzpicture}[anchorbase]
\draw[usual] (0,0) to (1,1);
\draw[usual] (1,0) to (0,1);
\end{tikzpicture}
,\quad
\begin{tikzpicture}[tqft/cobordism/.style={draw},tqft/every lower boundary component/.style={draw=black},anchorbase]
\pic[tqft/cylinder,genus=1];
\end{tikzpicture}
\leftrightsquigarrow
\begin{tikzpicture}[anchorbase]
\draw[usual,dotg] (0,0) to (0,1);
\end{tikzpicture}
,\quad
\text{\etc}
\end{gather*}
Such categories depend on the choice of a generating function as a quotient of two polynomials $p/q$, but for us only $k=\max\{\deg p+1,\deg q\}$ plays a role. More precisely, the coefficients $a_{i}$ of the Taylor expansion of $p/q$ are used to evaluate closed surfaces:
\begin{gather}\label{Eq:Evaluate}
\begin{tikzpicture}[anchorbase]
\draw[usual,dot=0] (0,0) to (0,0.1);
\draw[usual] (0,0.1) to (0,0.9);
\draw[usual,dot=1] (0,0.9) to (0,1);
\end{tikzpicture}
=a_{0},\quad
\begin{tikzpicture}[anchorbase]
\draw[usual,dot=0] (0,0) to (0,0.1);
\draw[usual,dotg] (0,0.1) to (0,0.9);
\draw[usual,dot=1] (0,0.9) to (0,1);
\end{tikzpicture}
=a_{1},\quad
\begin{tikzpicture}[anchorbase]
\draw[usual,dot=0] (0,0) to (0,0.1);
\draw[usual,dotg=0.25,dotg=0.75] (0,0.1) to (0,0.9);
\draw[usual,dot=1] (0,0.9) to (0,1);
\end{tikzpicture}
=a_{2},\quad
\dots,
\end{gather}
and $k$ is the degree of a minimal polynomial of the handle.
A special case for \changed{$p=t,q=1-x$}, so $k=1$, is the \emph{partition category} $\pamon[t]$, where all handles disappear
\begin{gather*}
\begin{tikzpicture}[anchorbase]
\draw[usual,dotg] (0,0) to (0,1);
\end{tikzpicture}
=
t\cdot
\begin{tikzpicture}[anchorbase]
\draw[usual] (0,0) to (0,1);
\end{tikzpicture}
,
\end{gather*}
and $t=a_{0}=a_{1}=\dots\in\C$ (generic in this paper) is the value of floating components, which was the category studied in \cite{De-cat-st}. From this we also get subcategories of the partition category that are related to the classical \emph{diagram monoids}, defining their diagrams. These monoids are the endomorphism monoids of the respective categories for $t=1$ (the list is taken from \cite{KhSiTu-monoidal-cryptography}\changed{; the subscript denotes the parameter, while the number in brackets indicates the object $n$, meaning that one takes some endomorphisms of $n\in\N$ circles $\bullet$}):
\begin{enumerate}[label=$\bullet$]

\item The \emph{partition monoid} \changed{$\pamon[1](n)$} of all diagrams of partitions of a $2n$-element set.

\item The \emph{rook-Brauer monoid} \changed{$\robrmon[1](n)$} consisting of all diagrams with components of size $1,2$.

\item The \emph{Brauer monoid} \changed{$\brmon[1](n)$} consisting of all diagrams with components of size $2$.

\item The \emph{rook monoid} \changed{$\romon[1](n)$} consisting of all diagrams with components of size $1,2$, and all partitions have at most one component
at the bottom and at most one at the top.

\item The \emph{symmetric group} \changed{$\sym[1](n)$} consisting of all matchings with components of size $1$.

\item \emph{Planar} versions of these: \changed{$\ppamon[1](n)$, $\probrmon[1](n)=\momon[1](n)$, $\pbrmon[1](n)=\tlmon[1](n)$, $\promon[1](n)$ and $\psym[1](n)\cong 1$} (the latter denotes the
trivial monoid). The planar rook-Brauer monoid is also called \emph{Motzkin monoid}, the planar Brauer monoid is also known as \emph{Temperley--Lieb monoid}, and the planar symmetric group is trivial.

\item Additionally, there are oriented versions of these. The one we will use is the \emph{oriented Brauer monoid} \changed{$\obrmon[1](n)$, which is $\brmon[1](n)$} with orientations.

\end{enumerate}

We also consider diagram categories that do not come from monoids. Precisely, the diagram categories that interpolate the categories of \ochanged{finite dimensional} complex representations \ochanged{in} $\repgl$, as studied in \cite{Kn-semisimple-tensor-gltfq} 
(these are diagram categories by \cite{EAH}, namely, these have partition diagrams plus extra generators).

\begin{Remark}
When referred to as \emph{interpolation categories} (potentially after taking some envelope), 
$\pamon[t]$, $\brmon[t]$ and $\obrmon[t]$ are often called 
Deligne(--Jones--Martin) categories \changed{$\repsnn[t]$}, $\reposp[t]$ (or $\repospt[t]$) and $\repgll[t]$, respectively, sometimes \changed{denoted} with an underline (similarly for \changed{$\repgl[t]$}). Since we focus on the semisimple case, there will be no difference between the additive or abelian versions.
\end{Remark}

Summarized with $\obstuff{X}$ denoting the chosen object and $\sim$ meaning \emph{asymptotically equal}, our results are \changed{outlined} in \autoref{Eq:Main}.

Following \cite{KhSiTu-monoidal-cryptography}, we establish these results using the Green relations and the cell structure of the associated monoids and diagram algebras, with the key ingredient being their \emph{sandwich cellular structure} \cite{Br-gen-matrix-algebras,TuVa-handlebody,Tu-sandwich-cellular}. This (new) perspective simplifies the problem significantly: in most cases, the first step of counting has already been carried out in semigroup theory \ochanged{(albeit from a very different viewpoint)} under the framework of Gelfand models \cite{HaRe-gelfand-diagrams}. The only exceptions are $\cob$ and $\repgl[t]$, which we analyze in detail. This approach yields exact formulas for $b_n$.  

The remaining task\ochanged{, deriving asymptotics from these exact formulas,} is nontrivial. However, as we show below, we have formalized the process so that it is computer-verified. For brevity, we have outsourced the computational aspects to \cite{GrTu}, where they are freely available (at least in 2025).  
\medskip

\noindent\textbf{Acknowledgments.}
We thank Johannes Flake for helpful comments and discussions\ochanged{, and the referee for a careful reading of the document and many helpful suggestions.}
JG also acknowledges financial support through the DFG project 531430348 and expresses gratitude for the support and hospitality provided by the Sydney Mathematical Research Institute (SMRI).
DT was supported by the ARC Future Fellowship FT230100489 and is deeply committed to the belief in the importance of open access.

\section{Schur--Weyl duality, sandwich cellular algebras and growth problems}\label{S:Cell}

For some field $\K$, assume that one has an additive Krull--Schmidt monoidal $\K$-linear category $\catstuff{C}$ with \ochanged{finite dimensional} hom-spaces, and an object $\obstuff{X}\in\catstuff{C}$. For simplicity, assume that $\catstuff{C}$ is semisimple. A version of Schur--Weyl duality implies that ($\obstuff{Y}$ is a simple summand of $\obstuff{X}^{\otimes n}$ that appears with multiplicity $m>0$ if and only if the semisimple algebra
$A_{n}=\End_{\catstuff{C}}(\obstuff{X}^{\otimes n})$ has a simple representation $L_{\obstuff{Y}}$ of dimension $\dim_{\K}L_{\obstuff{Y}}=m$), and there is a bijection between such $\obstuff{Y}$ and $L_{\obstuff{Y}}$.
In particular, we have:

\begin{Lemma}\label{L:Endo}
In the above setting, $A_{n}$ is semisimple and
\begin{gather*}
b_{n}=\sum_{L}\dim_{\K}L,\quad
(\text{sum over simple $A_{n}$-representations $L$}).
\end{gather*}
\end{Lemma}

\begin{proof}
Directly from the above discussion, which, in turn, can be justified as 
in \cite[Section 4C]{AnStTu-cellular-tilting}. See also 
\cite{Er-tensor-dimensions-symmetric-group,MR1698856}.
\end{proof}

We will now assume familiarity with sandwich cellular algebras \cite{Br-gen-matrix-algebras,TuVa-handlebody,Tu-sandwich-cellular}, or at least with some variation of it, most notably, \cite{FiGr-canonical-cases-brauer}.

\begin{Lemma}\label{L:Endo2}
If $A_{n}$ is a semisimple involutive sandwich cellular algebra with apex set $\mathcal{P}^{ap}$, bottom sets $\mathcal{B}_{\lambda}$ and sandwiches algebras $\mathcal{H}_{\lambda}$, then
\begin{gather*}
b_{n}=
\sum_{\lambda\in\mathcal{P}^{ap}}
\big(\#\mathcal{B}_{\lambda}\cdot
\sum_{L}\dim_{\K}L\big),\quad
(\text{inner sum over simple $\mathcal{H}_{\lambda}$-representations $L$}).
\end{gather*}
\end{Lemma}

\begin{proof}
Directly from \autoref{L:Endo} and the standard theory of sandwich cellular algebras.
\end{proof}

We will use \autoref{L:Endo2} throughout.

\section{Examples}\label{S:Examples}

We now go through the list of example in \autoref{S:Intro}. We will discuss $\cob$ and $\repgl[t]$ carefully. All other cases are similar to the count in $\cob$, so we just list the needed results for them.

\subsection{Cobordisms}\label{SS:Cob}

In the notation of the introduction, let $\cob[\infty]$ be the monoidal category 
with $\otimes$-generating object $\bullet$ and generating
$\circ$-$\otimes$-generating morphisms 
\begin{gather*}
\text{multiplication: }
\begin{tikzpicture}[anchorbase]
\draw[usual] (0,0) to (0.5,0.5) to (1,0);
\draw[usual] (0.5,0.5) to (0.5,1);
\end{tikzpicture}
,\quad
\text{comultiplication: }
\begin{tikzpicture}[anchorbase,yscale=-1]
\draw[usual] (0,0) to (0.5,0.5) to (1,0);
\draw[usual] (0.5,0.5) to (0.5,1);
\end{tikzpicture}
,\quad
\text{unit: }
\begin{tikzpicture}[anchorbase,yscale=-1]
\draw[usual,dot] (0,0) to (0,0.5);
\end{tikzpicture}
,\quad
\text{counit: }
\begin{tikzpicture}[anchorbase]
\draw[usual,dot] (0,0) to (0,0.5);
\end{tikzpicture}
,\quad
\text{crossing:}
\begin{tikzpicture}[anchorbase]
\draw[usual] (0,0) to (1,1);
\draw[usual] (1,0) to (0,1);
\end{tikzpicture}
,
\end{gather*}
modulo the $\circ$-$\otimes$-ideal that makes the crossing a symmetry and $\bullet$ a symmetric Frobenius object with the structure maps matching the nomenclature. See \cite{Ko-tqfts} for details.

The \emph{diagrammatic antiinvolution} ${}^{\ast}$ flips a cobordism \changed{upside-down}.
We call a diagram a \emph{merge diagram} if it contains only multiplications, counits and a minimal number of crossings. A \emph{split diagram} is a ${}^{\ast}$-flipped merge diagram. A \emph{dotted permutation diagram} contains only dots (=handles) and crossings. Here are some examples (flipping the left illustration gives a split diagram):
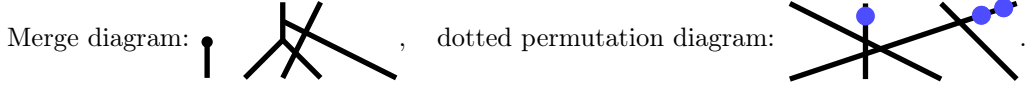
\begin{gather*}
\text{Merge diagram: }
\begin{tikzpicture}[anchorbase]
\draw[usual] (0,0) to (0.5,0.5) to (1,0);
\draw[usual] (0.5,0.5) to (0.5,1);
\draw[usual] (2,0) to (0.5,0.75);
\draw[usual] (0.5,0) to (1,1);
\draw[usual,dot] (-0.5,0) to (-0.5,0.5);
\end{tikzpicture}
,\quad
\text{dotted permutation diagram: }
\begin{tikzpicture}[anchorbase]
\draw[usual,dotg=0.85,dotg=0.95] (0,0) to (3,1);
\draw[usual,dotg=0.85] (1,0) to (1,1);
\draw[usual] (2,0) to (0,1);
\draw[usual] (3,0) to (2,1);
\end{tikzpicture}
.
\end{gather*}
Abusing notation, for a field $\K$ we denote the $\K$-linear extension of 
the cobordism category also by $\cob[\infty]$.
\ochanged{In recalling sandwich cellularity, we have the following picture for} \changed{$\cob[\infty]$}\ochanged{:}
\begin{gather}\label{Eq:CellsCob}
\begin{tikzpicture}[anchorbase,scale=1]
\draw[line width=0.75,color=black,fill=cream] (0,1) to (0.25,0.5) to (0.75,0.5) to (1,1) to (0,1);
\node at (0.5,0.75){$T$};
\draw[line width=0.75,color=black,fill=cream] (0.25,0) to (0.25,0.5) to (0.75,0.5) to (0.75,0) to (0.25,0);
\node at (0.51,0.25){$m$};
\draw[line width=0.75,color=black,fill=cream] (0,-0.5) to (0.25,0) to (0.75,0) to (1,-0.5) to (0,-0.5);
\node at (0.5,-0.25){$B$};
\end{tikzpicture}
\quad\text{where}\quad
\begin{aligned}
\begin{tikzpicture}[anchorbase,scale=1]
\draw[line width=0.75,color=black,fill=cream] (0,1) to (0.25,0.5) to (0.75,0.5) to (1,1) to (0,1);
\node at (0.5,0.75){$T$};
\end{tikzpicture}
&\text{ a split diagram,}
\\
\begin{tikzpicture}[anchorbase,scale=1]
\draw[white,ultra thin] (0,0) to (1,0);
\draw[line width=0.75,color=black,fill=cream] (0.25,0) to (0.25,0.5) to (0.75,0.5) to (0.75,0) to (0.25,0);
\node at (0.5,0.25){$m$};
\end{tikzpicture}
&\text{ a dotted permutation,}
\\
\begin{tikzpicture}[anchorbase,scale=1]
\draw[line width=0.75,color=black,fill=cream] (0,-0.5) to (0.25,0) to (0.75,0) to (1,-0.5) to (0,-0.5);
\node at (0.5,-0.25){$B$};
\end{tikzpicture}
&\text{ a merge diagram.}
\end{aligned}
\quad
\text{Involution: }
\left(\,
\begin{tikzpicture}[anchorbase,scale=1]
\draw[line width=0.75,color=black,fill=cream] (0,1) to (0.25,0.5) to (0.75,0.5) to (1,1) to (0,1);
\node at (0.5,0.75){$T$};
\draw[line width=0.75,color=black,fill=cream] (0.25,0) to (0.25,0.5) to (0.75,0.5) to (0.75,0) to (0.25,0);
\node at (0.51,0.25){$m$};
\draw[line width=0.75,color=black,fill=cream] (0,-0.5) to (0.25,0) to (0.75,0) to (1,-0.5) to (0,-0.5);
\node at (0.5,-0.25){$B$};
\end{tikzpicture}\,
\right)^{\ast}
=
\begin{tikzpicture}[anchorbase,scale=1]
\draw[line width=0.75,color=black,fill=cream] (0,1) to (0.25,0.5) to (0.75,0.5) to (1,1) to (0,1);
\node at (0.5,0.75){\changed{\reflectbox{\rotatebox{180}{$B$}}}};
\draw[line width=0.75,color=black,fill=cream] (0.25,0) to (0.25,0.5) to (0.75,0.5) to (0.75,0) to (0.25,0);
\node at (0.51,0.25){\reflectbox{\rotatebox{180}{$m$}}};
\draw[line width=0.75,color=black,fill=cream] (0,-0.5) to (0.25,0) to (0.75,0) to (1,-0.5) to (0,-0.5);
\node at (0.5,-0.25){\changed{\reflectbox{\rotatebox{180}{$T$}}}};
\end{tikzpicture}
.
\end{gather}
We call the existence of a spanning set with a \ochanged{decomposition} as above a precell structure and the respective algebras presandwich.

\begin{Lemma}\label{L:Cell}
The endomorphism algebras of $\cob[\infty]$ are involutive presandwich cellular with precell structure as in \autoref{Eq:CellsCob}.
\end{Lemma}

\begin{proof}
Immediate from \cite[Section 1.4.16]{Ko-tqfts}.
\end{proof}

Now fix two polynomials $p,q\in\K[x]$ \changed{with $q(0)\neq 0$}, and consider the Taylor expansion of $p(x)/q(x)=\sum_{i=0}^{\infty}a_{i}x^{i}$, and let $k=\max\{\deg p+1,\deg q\}$. Let 
$\cob[k]$ be the quotient of $\cob[\infty]$ by the $\circ$-$\otimes$-ideal generated by \autoref{Eq:Evaluate}. (The category $\cob[k]$ actually depends on $p,q$ but we suppress this in the notation.)

\begin{Example}\label{E:Partition}
\changed{For $p=t$ and $q=1-x$ we have $p(x)/q(x)=t(1+x+x^{2}+x^{3}+\dots)$ so that all closed surfaces in $\cob[k]$ evaluate to $t$.} We call the resulting category the \emph{partition category} \ochanged{$\pamon[t]$}.
\end{Example}

Let us denote the endomorphism monoid of $\bullet^{\otimes n}$ by $\cob[k](n)$. Recall the Ariki--Koike algebra (cyclotomic Hecke algebra) as defined in \cite{ArKo-hecke-algebra,BrMa-hecke,Ch-gelfandtzetlin}.

\begin{Lemma}\label{L:Cell2}
Dotted permutations in $\cob[k](n)$ span an algebra isomorphic to the Ariki--Koike algebra $A(n,k)$ on $n$ strands with a cyclotomic relation of degree $k$ and trivial quantum parameter. Moreover, \autoref{L:Cell} can be refined into a sandwich cell datum with $A(m,k)$ for $m\in\{0,\dots,n\}$ as the sandwiched algebras.
\end{Lemma}

\begin{proof}
This follows from \cite{KhOsKo-cobordisms}, {\eg} the text around (11) therein, which implies that the handles satisfy a minimal polynomial of degree $k$, and the same arguments as in \cite[Section 6]{TuVa-handlebody}.
\end{proof}

For the next statement we assume familiarity with the usual tableaux combinatorics, see, for example, \cite{DiJaMa-cyclotomic-q-schur,Ma-hecke-schur}.

\begin{Proposition}\label{P:Cob}
Let $\K$ be of characteristic $p$.
\begin{enumerate}

\item The set of apexes of $\K\cob[k](n)$ is $\{0,\dots,n\}$.

\item The \ochanged{finite dimensional} simple $\K\cob[k](n)$-modules of apex $m$, up to equivalence, are indexed by $p$-restricted $k$-multipartitions of $m$.

\item The dimensions of the cell representation for a $p$-restricted $k$-multipartition $\lambda$ of $m$ is
\begin{gather*}
\#\{\text{merge diagrams with $m$ top strands}\}
\cdot
\#\{\text{standard tableaux of shape $\lambda$}\}.
\end{gather*}

\end{enumerate}
Moreover, if $\K\cob[k](n)$ is semisimple, then the cell representations are simple.
\end{Proposition}

\begin{proof}
Immediate from the standard theory of sandwich cellular algebras as, {\eg}, in \cite{Tu-sandwich-cellular}, \autoref{L:Cell2} and 
the cell structure of $A(m,k)$ as, for example, in \cite[Theorem 3.26]{DiJaMa-cyclotomic-q-schur}.
\end{proof}

We now need some counting lemmas.

\begin{Lemma}\label{L:Count}
The number of merge diagrams from $n$ bottom strands to $m$ top strands is
\begin{gather*}
M_{n}^{m}=
\sum_{i=m}^{n}
\left\{\begin{matrix}n \\ i\end{matrix}\right\}
\binom{i}{m},
\end{gather*}
where the curly brackets denote the Stirling numbers of the second kind.
\end{Lemma}

\begin{proof}
A standard count that is independent of $k$, and therefore the same as in the partition category. Details are omitted; however, if the reader encounters difficulties, \cite[Section 4]{HaRe-gelfand-diagrams} provides helpful guidance.
\end{proof}

Let $\mathrm{STab}(m,k)$ is the set of standard $k$-multitableaux of $m$ and let $\#\mathrm{STab}(m,k)$ denote its size. Assume from now that $\K=\C$ and that $\cob[k]$ is semisimple. \changed{This holds, for example, if the $a_{i}$ are generic, {\eg} $a_{i}=t$ for a variable $t$ (giving the partition category for generic parameter).}

\begin{Lemma}
We have the formula
\begin{gather*}
b_{n}=\sum_{m=0}^{n}M_{n}^{m}\#\mathrm{STab}(m,k).
\end{gather*}
\end{Lemma}

\begin{proof}
By \autoref{P:Cob} and \autoref{L:Count}.
\end{proof}

\changed{Recall that the exponential generating function for a sequence $x_{n}$ is the generating function for $x_{n}/n!$.}

\begin{Lemma}\label{L:EGF}
\changed{The sequence} $b_{n}$ has exponential generating function $\exp(\frac{k}{2}\exp(2x)+\exp(x)-\frac{k+2}{2})$.
\end{Lemma}

\begin{proof}
We first observe that
\begin{gather*}
\#\mathrm{STab}(m,k)=k^{\lceil m/2\rceil}\sum_{i\in\N}Bes(m,i)k^{\lfloor m/2\rfloor-i},
\end{gather*}
where we use the Bessel numbers $Bes(m,i)=m!/(i!(n-2i)!2^{i})$. Thus, we get
\begin{gather*}
b_{n}=\sum_{m=0}^{n}\sum_{i\in\N}M_{n}^{m}k^{\lceil m/2\rceil}Bes(m,i)k^{\lfloor m/2\rfloor-i},
\end{gather*}
and we can use the same calculations as in \cite[Section 3]{calc} (which only uses the exponential generating functions for the Stirling and Bell numbers).
\end{proof}

\begin{Lemma}\label{L:Auto}
We have the asymptotic formula
\begin{gather*}
\changed{b_{n}\sim
\frac{\left(\frac{n}{z}\right)^{n+\frac{1}{2}}\exp\big(\frac{k}{2}\exp(2z+1) +\exp(z)-n-\frac{k+2}{2}\big)}{\sqrt{k\exp(2z)(2z+1)+\exp(z)(z+1))}}
,}
\\
\changed{z=\frac{W\left(\frac{2n}{k}\right)}{2} - \frac{1}{2k n^{\frac{1}{2}} \left( W\left(\frac{2n}{k}\right) + 1 \right) W\left(\frac{2n}{k}\right)^{-\frac{3}{2}} + \frac{2}{W\left(\frac{2n}{k}\right)} + 1},}
\end{gather*}
where $W$ is the Lambert W function.
\end{Lemma}

\begin{proof}
Having \autoref{L:EGF}, the proof of this is automatized, see for example \cite{GrTu,Ko}.
This works roughly as follows. Let $f$ be the exponential generating function. One then uses Hayman's method \cite{Hayman} and computes the asymptotic of $\lim_{x\to\infty}xf^{\prime}(x)/f(x)$.
\end{proof}

\begin{Theorem}
The formula in \autoref{Eq:Main} holds.
\end{Theorem}

\begin{proof}
The proof is also automatized, using \autoref{L:Auto} and the code 
on \cite{GrTu}. Essentially, Mathematica has a build in function for this purpose that does exact calculations. (We are referring to `Asymptotic'; Introduced in 2020 (12.1) | Updated in 2022 (13.2).) 
\end{proof}

\subsection{Partition algebras}

For \changed{$\pamon[t]$} we simply specialize $k=1$ in \autoref{SS:Cob}. 
See also \autoref{E:Partition}. The 
corresponding sequence for $b_{n}$ is \cite[A002872]{Oeis}. However, 
following \cite[A002872]{Oeis} one gets other, slightly nastier, formulas, namely:
\begin{gather*}
b_{n}\sim
\left( \frac{2n}{\text{W}(2n)} \right)^n \cdot \exp\left( \frac{n}{\text{W}(2n)} + \left( \frac{2n}{\text{W}(2n)} \right)^{\frac{1}{2}} - n - \frac{7}{4} \right) \Big/ \sqrt{1 + \text{W}(2n)},
\\
\sqrt[n]{b_{n}}\sim
\frac{2}{e}\cdot\frac{n}{\log 2n}\sqrt[\log 2n]{e}.
\end{gather*}
\changed{The term $\sqrt[\log 2n]{e}$ is created by Mathematica (and kept for that reason). It can be ignored.}

\subsection{Subalgebras of partition algebras}

For \ochanged{the current section and \autoref{SS:Oriented} we assume familiarity with the diagram monoids from \autoref{Eq:Main} and their respective categories.}
By an easy (and well known) diagrammatic argument we get \changed{that} \ochanged{$\ppamon[t](n)\cong\tlmon[t](2n)$}, we already discussed \changed{$\pamon[t]$, and $\psym[t]$} is trivial, so we do not need to address these cases.
Let us list the sandwich cellular bases for the remaining diagram categories:
\begin{gather*}
\momon,\tlmon,\promon\colon
\begin{aligned}
\begin{tikzpicture}[anchorbase,scale=1]
\draw[line width=0.75,color=black,fill=cream] (0,1) to (0.25,0.5) to (0.75,0.5) to (1,1) to (0,1);
\node at (0.5,0.75){$T$};
\end{tikzpicture}
&\text{ cup and unit diagrams,}
\\
\begin{tikzpicture}[anchorbase,scale=1]
\draw[white,ultra thin] (0,0) to (1,0);
\draw[line width=0.75,color=black,fill=cream] (0.25,0) to (0.25,0.5) to (0.75,0.5) to (0.75,0) to (0.25,0);
\node at (0.5,0.25){$m$};
\end{tikzpicture}
&\text{ an identity,}
\\
\begin{tikzpicture}[anchorbase,scale=1]
\draw[line width=0.75,color=black,fill=cream] (0,-0.5) to (0.25,0) to (0.75,0) to (1,-0.5) to (0,-0.5);
\node at (0.5,-0.25){$B$};
\end{tikzpicture}
&\text{ cap and counit diagram,}
\end{aligned}
\quad
\robrmon,\brmon,\romon\colon
\begin{aligned}
\begin{tikzpicture}[anchorbase,scale=1]
\draw[line width=0.75,color=black,fill=cream] (0,1) to (0.25,0.5) to (0.75,0.5) to (1,1) to (0,1);
\node at (0.5,0.75){$T$};
\end{tikzpicture}
&\text{ cup and unit diagrams,}
\\
\begin{tikzpicture}[anchorbase,scale=1]
\draw[white,ultra thin] (0,0) to (1,0);
\draw[line width=0.75,color=black,fill=cream] (0.25,0) to (0.25,0.5) to (0.75,0.5) to (0.75,0) to (0.25,0);
\node at (0.5,0.25){$m$};
\end{tikzpicture}
&\text{ a permutation,}
\\
\begin{tikzpicture}[anchorbase,scale=1]
\draw[line width=0.75,color=black,fill=cream] (0,-0.5) to (0.25,0) to (0.75,0) to (1,-0.5) to (0,-0.5);
\node at (0.5,-0.25){$B$};
\end{tikzpicture}
&\text{ cap and counit diagram.}
\end{aligned}
\end{gather*}
For the symmetric group $\sym$ the sandwich structure is trivial.

From this one gets explicit formulas, matching the ones in \cite[Section 4]{HaRe-gelfand-diagrams}.
Here is the list of \ochanged{the} remaining sequences:
\begin{gather*}
\momon\colon\text{\cite[A005773]{Oeis}},
\quad
\robrmon\colon\text{\cite[A000898]{Oeis}},
\\
\tlmon\colon\text{\cite[A000984]{Oeis}},
\quad
\brmon\colon\text{\cite[A047974]{Oeis}},
\\
\promon\colon\text{\cite[A000079]{Oeis}},
\quad
\romon\colon\text{\cite[A005425]{Oeis}},
\\
\sym\colon\text{\cite[A000085]{Oeis}}.
\end{gather*}
Let us just focus on $\brmon$; the others being similar. In this case the exponential generating function 
is $\exp(x^{2}+x)$, and Mathematica gives
\begin{gather*}
b_{n}\sim
2^{\frac{n}{2} - \frac{1}{2}} \exp\left(\sqrt{\frac{n}{2}} - \frac{n}{2} - \frac{1}{8}\right) n^{\frac{n}{2}},
\\
\sqrt[n]{b_{n}}\sim
\frac{\sqrt{2}}{\sqrt{e}}\cdot\sqrt{n}
.
\end{gather*}
This completes the proof.

\subsection{Oriented Brauer algebras}\label{SS:Oriented}

The left case in \autoref{Eq:Main} is the same as for $\sym$. For the right case, we observe that crossings provide isomorphisms, so we can reorder $(\uparrow\downarrow)^{n}$ to $n$ upwards pointing arrows followed by $n$ downwards pointing arrows. Thus, the sandwich structure is (or rather, can be arranged to be)
as follows.
\begin{gather*}
\obrmon\colon
\begin{aligned}
\begin{tikzpicture}[anchorbase,scale=1]
\draw[line width=0.75,color=black,fill=cream] (0,1) to (0.25,0.5) to (0.75,0.5) to (1,1) to (0,1);
\node at (0.5,0.75){$T$};
\end{tikzpicture}
&\text{ cup diagrams passing the middle,}
\\
\begin{tikzpicture}[anchorbase,scale=1]
\draw[white,ultra thin] (0,0) to (1,0);
\draw[line width=0.75,color=black,fill=cream] (0.25,0) to (0.25,0.5) to (0.75,0.5) to (0.75,0) to (0.25,0);
\node at (0.5,0.25){$m$};
\end{tikzpicture}
&\text{ a $\uparrow$ permutation and a $\downarrow$ permutation,}
\\
\begin{tikzpicture}[anchorbase,scale=1]
\draw[line width=0.75,color=black,fill=cream] (0,-0.5) to (0.25,0) to (0.75,0) to (1,-0.5) to (0,-0.5);
\node at (0.5,-0.25){$B$};
\end{tikzpicture}
&\text{ cap diagram passing the middle.}
\end{aligned}
\end{gather*}
The set of cap diagrams passing the middle can be identified with 
matchings of size $n-k$ of $\{1,\dots,n\}$ with $\{1^{\prime},\dots,n^{\prime}\}$, which has size $\binom{n}{k}^{2}(n-k)!$.
Moreover, let $cd(k)$ denote the sum of the dimensions of simple $\sym[k]$-representations, {\ie} the sequence \cite[A000085]{Oeis}, which has exponential generating function $\exp(\frac{1}{2}x^{2}+x)$. A calculation, using the recursion 
$cd(n)=cd(n-1)+(n-1)cd(n-2)$, shows
\begin{gather*}
cd(m+n)=\sum_{k\geq 0}k!\binom{m}{k}\binom{n}{k}cd(m-k)cd(n-k).
\end{gather*}
Taking everything together, we get
\begin{gather*}
b_{n}=\sum_{k=0}^{n}cd(k)^{2}
\binom{n}{k}^{2}(n-k)!=cd(2n),
\end{gather*}
and Mathematica proves 
\begin{gather*}
b_{n}\sim
n^{n}2^{n-1/2}
\exp(-n+\sqrt{2n}-1/4)(1 + 7/(24\sqrt{2n})),
\end{gather*}
and the result in \autoref{Eq:Main} itself.

\begin{Remark}
For the lover of diagrammatics, as an alternative argument one could reorder $(\uparrow\downarrow)^{n}$ to $n$ upwards pointing arrows followed by $n$ downwards pointing arrows, as above, and then use that caps and cups are also invertible operations to bend the diagrams to look like $\sym[2n]$.
\end{Remark}

\subsection{The general linear group over a finite field}

Since the case of $\Rep\big( \GL_t(\F_q) \big)$ is the least studied, we provide more details than in the other cases. A more thorough treatment can be found in \cite{GrTu}.

Throughout this section, we fix a prime power $q$.
We start by recalling the definition of the interpolation category $\Rep\big( \GL_t(\F_q) \big)$, following \cite{Knop_tensor_envelopes}, and by explaining the sandwich cellular structure for endomorphism algebras in $\Rep\big(\GL_t(\F_q)\big)$.

\subsubsection{The definition of \texorpdfstring{$\Rep\big( \GL_t(\F_q) \big)$}{Rep( GLt(Fq) )}}

For $m,n \in \Z_{\geq 0}$, let us write $\mathrm{Gr}(\F_q^m,\F_q^n)$ for the set of linear spaces of $\F_q^m \oplus \F_q^n$.
For two subspaces $V \in \mathrm{Gr}(\F_q^m,\F_q^n)$ and $W \in \mathrm{Gr}(\F_q^\ell,\F_q^m)$, the convolution $V \star W \in \mathrm{Gr}(\F_q^\ell,\F_q^n)$ is defined as the image in $\F_q^\ell \oplus \F_q^n$ of the pullback $W \times_{\F_q^m} V$, as in the following diagram \changed{(using the standard projections)}.
\begin{center}
\begin{tikzcd}
    & & W \times_{\F_q^m} V \ar[ld] \ar[rd] \arrow[dashed,two heads]{d}{e} & & \\
    & W \ar[ld] \ar[rd] & V \star W \ar[lld,dashed] \ar[rrd,dashed] & V \ar[ld] \ar[rd] & \\
    \F_q^\ell & & \F_q^m & & \F_q^n
\end{tikzcd}
.
\end{center}
Furthermore, we write $d(V,W) = \dim \ker(e)$ for the dimension of the kernel of the canonical epimorphism $e \colon W \times_{\F_q^m} V \to V \star W$.

\begin{Definition}
    For a complex parameter $t \in \C$, the $\C$-linear category $\Rep\big( \GL_t(\F_q) \big)_0$ has objects $\Z_{\geq 0}$ and
    %$\Hom$-spaces
    \[ \Hom_{\GL_t(\F_q)}(m,n) = \C \mathrm{Gr}(\F_q^m,\F_q^n), \]
    for $m,n \in \Z_{\geq 0}$, the $\C$-vector space with basis $\mathrm{Gr}(\F_q^m,\F_q^n)$.
    The composition of homomorphisms is defined via
    \[ V \circ W = t^{d(V,W)} \cdot V \star W, \]
    for $V \in \mathrm{Gr}(\F_q^m,\F_q^n)$ and $W \in \mathrm{Gr}(\F_q^\ell,\F_q^m)$, 
    extended by \ochanged{bilinearly}.
    
    The category $\Rep\big( \GL_t(\F_q) \big)$ is defined as the additive Karoubi envelope of $\Rep\big(\GL_t(\F_q)\big)_0$ (i.e.\ the idempotent completion of the additive envelope).
\end{Definition}

The category $\Rep\big( \GL_t(\F_q) \big)$ has a \ochanged{(canonical)} $\C$-linear symmetric monoidal structure, with tensor unit $\mathbf{1} = 0$ and tensor product given by $m \otimes n = m+n$ and
\[ V \otimes W = V \oplus W \subseteq ( \F_q^m \oplus \F_q^n) \oplus ( \F_q^{m'} \oplus \F_q^{n'} ), \]
for $V \in \mathrm{Gr}(\F_q^m,\F_q^{m'})$ and $W \in \mathrm{Gr}(\F_q^n,\F_q^{n'})$, where the tensor product on the \ochanged{left-hand} side is taken in $\Rep\big(\GL_t(\F_q)\big)$ and the direct sum on the \ochanged{right-hand} side is taken in $\F_q$-vector spaces.
Equipped with this monoidal structure, the category $\Rep\big(\GL_t(\F_q)\big)$ is rigid, with evaluation and coevaluation maps both given by the diagonal subspace $\mathrm{diag}(m) = \{ (v,v) \mid v \in \F_q^m \} \subseteq \F_q^{2m}$, viewed either as a homomorphism $m \otimes m \to 0$ or as a homomorphism $0 \to m \otimes m$.
The dual of $V \in \mathrm{Gr}(\F_q^m,\F_q^n)$, viewed as a homomorphism $m \to n$ in $\Rep\big(\GL_t(\F_q)\big)$, is the subspace 
\[ V^* = \{ (v,w) \in \F_q^n \oplus \F_q^m \mid (w,v) \in V \} \in \mathrm{Gr}(\F_q^n,\F_q^m) , \]
viewed as a homomorphism $n \to m$.

\begin{Remark}
    For an $n \times m$-matrix $M$ with entries in $\F_q$, the graph $G(M) = \{ ( v, Mv ) \mid v \in \F_q^m \}$ is a subspace of $\F_q^m \oplus \F_q^n$, and for an $m \times \ell$-matrix $M'$ with entries in $\F_q$, one easily checks that $G(M) \circ G(M') = G(M M')$.
    %For every $\F_q$-linear map $f \colon \F_q^m \to \F_q^n$, the graph $G(f) = \{ ( v, f(v) ) \mid v \in \F_q^m \}$ is a subspace of $\F_q^m \oplus \F_q^n$, and for another $\F_q$-linear map $f' \colon \F_q^\ell \to \F_q^m$, one easily checks that $G(f) \circ G(f') = G(f \circ f')$.
    %(In other words, $G$ defines a functor from the category of \ochanged{finite dimensional} $\F_q$-vector spaces to $\Rep\big( \GL_t(\F_q) \big)$.)
    In particular, the group algebra $\C[\GL_n(\F_q)]$ is a subalgebra of the endomorphism ring $\End_{\GL_t(\F_q)}(n)$.
\end{Remark}

\subsubsection{Sandwich cellular structure}
\label{subsubsec:RepGLtFqsandwichcellular}

The sandwich cellular structure of endomorphism algebras in $\Rep\big(\GL_t(\F_q)\big)$ relies on so-called \emph{core factorizations}, as explained (in a more general setting) in \cite[Section 5]{Knop_tensor_envelopes}.
For $m \geq n \geq 0$, let us write $\mathrm{Gr}^\mathrm{IS}(\F_q^m,\F_q^n)$ for the set of subspaces $V$ of $\F_q^m \oplus \F_q^n$ such that the projection $V \to \F_q^m$ is injective and the projection $V \to \F_q^n$ is surjective, and similarly define \changed{$\mathrm{Gr}^\mathrm{SI}(\F_q^n,\F_q^m)$} as the set of subspaces $V$ of $\F_q^n \oplus \F_q^m$ such that the projection $V \to \F_q^n$ is surjective and the projection $V \to \F_q^m$ is injective.
(The notation IS stands for ``injective-surjective'', and SI is for ``surjective-injective''.)

\begin{Lemma}[Core factorization]
    Let $m,n \in \Z_{\geq 0}$ and $V \in \mathrm{Gr}(\F_q^m,\F_q^n)$.
    \begin{enumerate}
    \item There is some $k \leq \min\{ m , n \}$ and $V_1 \in \mathrm{Gr}^\mathrm{IS}(\F_q^m,\F_q^k)$ and $V_2 \in \mathrm{Gr}^\mathrm{SI}(\F_q^k,\F_q^n)$ such that $V = V_2 \circ V_1$.
    \item For any $\ell \leq \min\{m,n\}$ and $V_1' \in \mathrm{Gr}^\mathrm{IS}(\F_q^m,\F_q^\ell)$ and $V_2' \in \mathrm{Gr}^\mathrm{SI}(\F_q^\ell,\F_q^n)$ such that $V = V_2' \circ V_1'$, we have $\ell = k$ and there is $M \in \GL_k(\F_q)$ such that $V_1' = G(M) \circ V_1$ and $V_2 = V_2' \circ G(M)$.
    \end{enumerate}
\end{Lemma}

In other words, the composition of homomorphisms in $\Rep\big(\GL_t(\F_q)\big)$ gives rise to a bijection
\[ \bigsqcup_{k \leq \min\{m,n\} } \big( \mathrm{Gr}^\mathrm{IS}(\F_q^m,\F_q^k) \times \mathrm{Gr}^\mathrm{SI}(\F_q^k,\F_q^n) \big) \, \big/ \, \GL_k(\F_q) ~ \xrightarrow{ ~1:1~ } ~ \mathrm{Gr}(\F_q^m,\F_q^n) , \]
where $M \in \GL_k(\F_q)$ acts on $\mathrm{Gr}^\mathrm{IS}(\F_q^m,\F_q^k) \times \mathrm{Gr}^\mathrm{SI}(\F_q^k,\F_q^n)$ via $(V_1,V_2) \mapsto ( G(M) \circ V_1 , V_2 \circ G(M^{-1}) )$.
Further note that $\GL_k(\F_q)$ acts freely on $\mathrm{Gr}^\mathrm{IS}(\F_q^m,\F_q^k)$ and $\mathrm{Gr}^\mathrm{SI}(\F_q^k,\F_q^n)$ via $(M,V_1) \mapsto G(M) \circ V_1$ and $(M,V_2) \mapsto V_2 \circ G(M^{-1})$, respectively.
In particular, for any fixed sets of representatives $\mathrm{Gr}^\mathrm{IS}_0(\F_q^m,\F_q^k)$ and $\mathrm{Gr}^\mathrm{SI}_0(\F_q^k,\F_q^n)$ for the $\GL_k(\F_q)$-orbits in $\mathrm{Gr}^\mathrm{IS}(\F_q^m,\F_q^k)$ and $\mathrm{Gr}^\mathrm{SI}(\F_q^k,\F_q^n)$, respectively, we have a bijection
\[ \mathrm{Gr}^\mathrm{IS}_0(\F_q^m,\F_q^k) \times \GL_k(\F_q) \times \mathrm{Gr}^\mathrm{SI}_0(\F_q^k,\F_q^n) ~ \xrightarrow{~ 1:1 ~} ~ \big( \mathrm{Gr}^\mathrm{IS}(\F_q^m,\F_q^k) \times \mathrm{Gr}^\mathrm{SI}(\F_q^k,\F_q^n) \big) \, \big/ \, \GL_k(\F_q) , \]
which sends a triple $(V_1,M,V_2)$ to the $\GL_k(\F_q)$-orbit of $( G(M) \circ V_1 , V_2 )$.
These observations essentially supply all of the data that makes up a sandwich cellular structure on the endomorphism algebras in $\Rep\big(\GL_t(\F_q)\big)$, following \cite[Definition 2A.3]{Tu-sandwich-cellular}.
Namely, for $n \in \Z_{\geq 0}$, the sandwich cell datum $( \mathcal{P} , (\mathcal{T},\mathcal{B}) , (\mathcal{H}_k,B_k) , C )$ for the algebra $\mathcal{A} = \End_{\GL_t(\F_q)}(n)$ is defined as follows:
\begin{itemize}
    \item The middle poset is $\mathcal{P} = \{ k \in \Z_{\geq 0} \mid k \leq n \}$, endowed with the usual partial \ochanged{(in this case even total)} order on integers.
    \item The top and bottom sets are given by $\mathcal{T} = \bigsqcup_{k \leq n} \mathcal{T}(k)$ and $\mathcal{B} = \bigsqcup_{k \leq n} \mathcal{B}(k)$, with
    \[ \mathcal{T}(k) = \mathrm{Gr}_0^\mathrm{SI}(\F_q^k,\F_q^n) , \hspace{2cm} \mathcal{B}(k) = \mathrm{Gr}_0^\mathrm{IS}(\F_q^n,\F_q^k) . \]
    \item The sandwiched algebras are the group algebras $\mathcal{H}_k = \C[\GL_k(\F_q)]$, with a fixed choice of basis given by $B_k = \{ M \mid M \in \GL_k(\F_q) \}$.
    \item The map $C \colon \bigsqcup_{k \leq n} \mathcal{T}(k) \times B_k \times \mathcal{B}(k) \longrightarrow \mathcal{A}$ indexing the sandwich cellular basis is given by
    \[ (V_1,M,V_2) \longmapsto c_{V_1,V_2}^M = V_1 \circ G(M) \circ V_2 . \]
\end{itemize}
It is clear from the above discussion that $C$ indexes a basis of $\mathcal{A} = \End_{\GL_t(\F_q)}(n)$.
The cell modules from ($\mathrm{AC}_2$) in \cite{Tu-sandwich-cellular} are the $\mathcal{A} \text{-} \mathcal{H}_k$-bimodule $\Delta(k) = \C \mathrm{Gr}^{\mathrm{SI}}(\F_q^n,\F_q^k)$ and the $\mathcal{H}_k \text{-} \mathcal{A}$-bimodule $\nabla(k) = \C \mathrm{Gr}^{\mathrm{IS}}(\F_q^k,\F_q^n)$.
Checking the remaining axioms in \cite[Definition 2A.3]{Tu-sandwich-cellular} is straightforward, but somewhat tedious, and we leave the details to the reader.

Summarized:
\begin{gather*}
\Rep\big(\GL_t(\F_q)\big)\colon
\begin{aligned}
\begin{tikzpicture}[anchorbase,scale=1]
\draw[line width=0.75,color=black,fill=cream] (0,1) to (0.25,0.5) to (0.75,0.5) to (1,1) to (0,1);
\node at (0.5,0.75){$T$};
\end{tikzpicture}
&\text{ SI subspaces $\mathrm{Gr}_0^\mathrm{SI}(\F_q^k,\F_q^n)$,}
\\
\begin{tikzpicture}[anchorbase,scale=1]
\draw[white,ultra thin] (0,0) to (1,0);
\draw[line width=0.75,color=black,fill=cream] (0.25,0) to (0.25,0.5) to (0.75,0.5) to (0.75,0) to (0.25,0);
\node at (0.5,0.25){$m$};
\end{tikzpicture}
&\text{ the group algebra $\C[\GL_k(\F_q)]$,}
\\
\begin{tikzpicture}[anchorbase,scale=1]
\draw[line width=0.75,color=black,fill=cream] (0,-0.5) to (0.25,0) to (0.75,0) to (1,-0.5) to (0,-0.5);
\node at (0.5,-0.25){$B$};
\end{tikzpicture}
&\text{ IS subspaces $\mathrm{Gr}_0^\mathrm{IS}(\F_q^n,\F_q^k)$.}
\end{aligned}
\end{gather*}

\subsubsection{Counting direct summands}

From now on, we assume that $q$ is odd \changed{(we need this assumption to compute the sequence $cd(k)$ defined below)}, and that $t \in \C \setminus \{ q^k \mid k \in \Z_{\geq 0} \}$, so that the category $\Rep\big(\GL_t(\F_q)\big)$ is semisimple by \cite[Theorem 8.8 and Example 4 in Section 8]{Knop_tensor_envelopes}.
As before, we consider the sequence
\[ b_n =\#\text{indecomposable summands in $1^{\otimes n}$ counted with multiplicities} . \]
By \autoref{L:Endo2} and the discussion in \autoref{subsubsec:RepGLtFqsandwichcellular}, we have
    \[ b_n = \sum_{k \leq n} \Big( \# \mathcal{B}(k) \cdot
\sum_{\chi} \chi(1) \Big) , \]
where $\mathcal{B}(k) = \mathrm{Gr}^\mathrm{IS}_0(\F_q^n,\F_q^k)$ and $\chi$ runs over the \ochanged{simple} characters of $\GL_k(\F_q)$.
Therefore, in order to explicitly determine $b_n$, we need to\ochanged{:}
\begin{enumerate}
    \item compute the sum $cd(k) = \sum_\chi \chi(1)$ of the degrees of the \ochanged{simple} characters of $\GL_k(\F_q)$;
    \item count the number of elements in $\mathrm{Gr}^\mathrm{IS}_0(\F_q^n,\F_q^k) \cong \mathrm{Gr}^\mathrm{IS}(\F_q^n,\F_q^k) \, / \, \GL_k(\F_q)$.
\end{enumerate}
Since we assume that $q$ is odd, $cd(k)$ equals the number of symmetric matrices in $\GL_k(\F_q)$ by \cite[Theorem 3]{Gow_finite_general_linear} \changed{(here $q$ is required to be odd)}, and so by \cite[Theorem 2]{MacWilliams_orthogonal_finite_field}, we have
\[ cd(k) = \# \text{symmetric matrices in } \GL_k(\F_q) = q^{\binom{k+1}{2}} \cdot \prod_{ \substack{1\leq i \leq k\\i \text{ odd}} } ( 1 - \tfrac{1}{q^i} ) . \]
Furthermore, arguing as in \autoref{subsubsec:RepGLtFqsandwichcellular}, one sees that there is a bijection
\[ \mathrm{Gr}^\mathrm{IS}(\F_q^n,\F_q^k) \, / \, \GL_k(\F_q) ~ \xrightarrow{~1:1~} ~ \bigsqcup_{k \leq \ell \leq n} \mathrm{Gr}_\ell(\F_q^n) \times \mathrm{Gr}_k(\F_q^\ell) , \]
where $\mathrm{Gr}_b(\F_q^a)$ denotes the set of $b$-dimensional subspaces in $\F_q^a$.
It is well-known that the number of elements in $\mathrm{Gr}_b(\F_q^a)$ is the (Gaussian) $q$-binomial coefficient
\[ \# \mathrm{Gr}_b(\F_q^a) = \left[ \begin{smallmatrix} a \\ b \end{smallmatrix} \right]_q = \frac{[a]_q [a-1]_q \cdots [a-b+1]_q}{[b]_q [b-1]_q \cdots [1]_q} , \]
with $[j]_q = \tfrac{q^j-1}{q-1}$ for $j \in \Z_{\geq 0}$.
In conclusion, we have
\begin{align*}
    b_n & = \sum_{k \leq n} \# \mathrm{Gr}^\mathrm{IS}_0(\F_q^n,\F_q^k) \cdot cd(k) \\
    & = \sum_{k \leq n} \Big( \sum_{k \leq \ell \leq n} \left[ \begin{smallmatrix} n \\ \ell \end{smallmatrix} \right]_q \left[ \begin{smallmatrix} \ell \\ k \end{smallmatrix} \right]_q \Big) \cdot q^{\binom{k+1}{2}} \cdot \prod_{\substack{1 \leq i \leq k\\i \text{ odd}}} ( 1 - \tfrac{1}{q^i} ) \\
    & = \prod_{k=1}^n (q^k+1) .
\end{align*}
\changed{The last equality can be computer verified or explicitly proven by hand (although the computation is rather nasty); see \cite{GrTu} for more details, {\eg} Proposition 3E.6 in the additional details file or the standalone \TeX~document on that website.}
If we set
\[ c = \lim_{n \to \infty} \prod_{k=1}^n (1 + \tfrac{1}{q^k}) = \tfrac{1}{2}\mathrm{QPochhammer}(-1,1/q)_{\infty} \]
(the q-Pochhammer symbol) then it follows that
\[ b_n \sim c \cdot q^{\frac{n(n+1)}{2}} , \hspace{2cm} \sqrt[n]{b_{n}} \sim q^{\frac{n+1}{2}} , \]
and this establishes the claim in \eqref{Eq:Main}.


\begin{thebibliography}{CEOT24}

\bibitem[AST18]{AnStTu-cellular-tilting}
H.H.~Andersen, C.~Stroppel, and D.~Tubbenhauer.
\newblock Cellular structures using {$\mathrm{U}_q$}-tilting modules.
\newblock {\em Pacific J. Math.}, 292(1):21--59, 2018.
\newblock URL: \url{https://arxiv.org/abs/1503.00224}, \href
  {https://doi.org/10.2140/pjm.2018.292.21}
  {\path{doi:10.2140/pjm.2018.292.21}}.

\bibitem[AK94]{ArKo-hecke-algebra}
S.~Ariki and K.~Koike.
\newblock A {H}ecke algebra of {$(\mathbb{Z}/r\mathbb{Z})\wr\mathfrak{S}_{n}$}
  and construction of its irreducible representations.
\newblock {\em Adv. Math.}, 106(2):216--243, 1994.
\newblock \href {https://doi.org/10.1006/aima.1994.1057}
  {\path{doi:10.1006/aima.1994.1057}}.

\bibitem[BM93]{BrMa-hecke}
\changed{M.~Brou{\'e} and G.~Malle.}
\newblock Zyklotomische {H}eckealgebren.
\newblock Number 212, pages 119--189. 1993.
\newblock Repr{\'e}sentations unipotentes g{\'e}n{\'e}riques et blocs des
  groupes r{\'e}ductifs finis.

\bibitem[Bro55]{Br-gen-matrix-algebras}
\changed{W.P.~Brown.}
\newblock Generalized matrix algebras.
\newblock {\em Canadian J. Math.}, 7:188--190, 1955.
\newblock \href {https://doi.org/10.4153/CJM-1955-023-2}
{\path{doi:10.4153/CJM-1955-023-2}}.

\bibitem[Che87]{Ch-gelfandtzetlin}
I.V.~Cherednik.
\newblock A new interpretation of {G}elfand--{T}zetlin bases.
\newblock {\em Duke Math. J.}, 54(2):563--577, 1987.
\newblock \href {https://doi.org/10.1215/S0012-7094-87-05423-8}
  {\path{doi:10.1215/S0012-7094-87-05423-8}}.

\bibitem[CEOT24]{CoEtOsTu-growth-fractal}
K.~Coulembier, P.~Etingof, V.~Ostrik, and D.~Tubbenhauer.
\newblock Fractal behavior of tensor powers of the two dimensional space in
  prime characteristic.
\newblock 2024.
\newblock To appear in Contemp. Math.
\newblock URL: \url{https://arxiv.org/abs/2405.16786}.

\bibitem[COT24]{CoOsTu-growth}
K.~Coulembier, V.~Ostrik, and D.~Tubbenhauer.
\newblock Growth rates of the number of indecomposable summands in tensor
  powers.
\newblock {\em Algebr. Represent. Theory}, 27(2):1033--1062, 2024.
\newblock URL: \url{https://arxiv.org/abs/2301.00885}, \href
  {https://doi.org/10.1007/s10468-023-10245-7}
  {\path{doi:10.1007/s10468-023-10245-7}}.

\bibitem[Del07]{De-cat-st}
P.~Deligne.
\newblock La cat{\'e}gorie des repr{\'e}sentations du groupe sym{\'e}trique
  {$S_t$}, lorsque {$t$} n'est pas un entier naturel.
\newblock In {\em Algebraic groups and homogeneous spaces}, volume~19 of {\em
  Tata Inst. Fund. Res. Stud. Math.}, pages 209--273. Tata Inst. Fund. Res.,
  Mumbai, 2007.

\bibitem[DJM98]{DiJaMa-cyclotomic-q-schur}
R.~Dipper, G.~James, and A.~Mathas.
\newblock Cyclotomic {$q$}-{S}chur algebras.
\newblock {\em Math. Z.}, 229(3):385--416, 1998.
\newblock \href {https://doi.org/10.1007/PL00004665}
  {\path{doi:10.1007/PL00004665}}.

\bibitem[EAH25]{EAH}
I.~Entova-Aizenbud and T.~Heidersdorf.
\newblock Deligne categories and representations of the finite general linear
  group, part {1}: universal property.
\newblock \ochanged{{\em Transform. Groups}, Volume 30, pages 633--698, (2025).}
\newblock URL: \url{https://arxiv.org/abs/2208.00241}, \href
  {https://doi.org/10.1007/s00031-023-09840-1}
  {\path{doi:10.1007/s00031-023-09840-1}}.

\bibitem[Erd95]{Er-tensor-dimensions-symmetric-group}
K.~Erdmann.
\newblock Tensor products and dimensions of simple modules for symmetric
  groups.
\newblock {\em Manuscripta Math.}, 88(3):357--386, 1995.
\newblock \href {https://doi.org/10.1007/BF02567828}
  {\path{doi:10.1007/BF02567828}}.

\bibitem[FG95]{FiGr-canonical-cases-brauer}
S.~Fishel and I.~Grojnowski.
\newblock Canonical bases for the {B}rauer centralizer algebra.
\newblock {\em Math. Res. Lett.}, 2(1):15--26, 1995.
\newblock \href {https://doi.org/10.4310/MRL.1995.v2.n1.a3}
  {\path{doi:10.4310/MRL.1995.v2.n1.a3}}.

\bibitem[Gow83]{Gow_finite_general_linear}
R.~Gow.
\newblock Properties of the characters of the finite general linear group
  related to the transpose-inverse involution.
\newblock {\em Proc. Lond. Math. Soc. (3)}, 47:493--506, 1983.
\newblock \href {https://doi.org/10.1112/plms/s3-47.3.493}
  {\path{doi:10.1112/plms/s3-47.3.493}}.

\bibitem[GT25]{GrTu}
J.~Gruber and D.~Tubbenhauer.
\newblock {C}ode and more for the paper {G}rowth problems in diagram
  categories.
\newblock 2025.
\newblock \url{https://github.com/dtubbenhauer/growth-diagram}.

\bibitem[HR15]{HaRe-gelfand-diagrams}
T.~Halverson and M.~Reeks.
\newblock Gelfand models for diagram algebras.
\newblock {\em J. Algebraic Combin.}, 41(2):229--255, 2015.
\newblock URL: \url{https://arxiv.org/abs/1302.6150}, \href
  {https://doi.org/10.1007/s10801-014-0534-5}
  {\path{doi:10.1007/s10801-014-0534-5}}.

\bibitem[Hay56]{Hayman}
W.K.~Hayman.
\newblock A generalisation of {S}tirling's formula.
\newblock {\em J. Reine Angew. Math.}, 196:67--95, 1956.
\newblock \href {https://doi.org/10.1515/crll.1956.196.67}
  {\path{doi:10.1515/crll.1956.196.67}}.

\bibitem[He24]{He-1}
D.~He.
\newblock Growth problems for representations of finite groups.
\newblock 2024.
\newblock To appear in Ark. Mat.
\newblock URL: \url{https://arxiv.org/abs/2408.04196}.

\bibitem[HT25]{He-2}
D.~He and D.~Tubbenhauer.
\newblock Growth problems for representations of finite monoids.
\newblock 2025.
\newblock URL: \url{https://arxiv.org/abs/2502.02849}.

\bibitem[KOK22]{KhOsKo-cobordisms}
M.~Khovanov, V.~Ostrik, and Y.~Kononov.
\newblock Two-dimensional topological theories, rational functions and their
  tensor envelopes.
\newblock {\em Selecta Math. (N.S.)}, 28(4):Paper No. 71, 68, 2022.
\newblock URL: \url{https://arxiv.org/abs/2011.14758}, \href
  {https://doi.org/10.1007/s00029-022-00785-z}
  {\path{doi:10.1007/s00029-022-00785-z}}.

\bibitem[KS24]{KhSa-cobordisms}
M.~Khovanov and R.~Sazdanovic.
\newblock Bilinear pairings on two-dimensional cobordisms and generalizations
  of the {D}eligne category.
\newblock {\em Fund. Math.}, 264(1):1--20, 2024.
\newblock URL: \url{https://arxiv.org/abs/2007.11640}, \href
  {https://doi.org/10.4064/fm283-8-2023} {\path{doi:10.4064/fm283-8-2023}}.

\bibitem[KST24]{KhSiTu-monoidal-cryptography}
M.~Khovanov, M.~Sitaraman, and D.~Tubbenhauer.
\newblock Monoidal categories, representation gap and cryptography.
\newblock {\em Trans. Amer. Math. Soc. Ser. B}, 11:329--395, 2024.
\newblock URL: \url{https://arxiv.org/abs/2201.01805}, \href
  {https://doi.org/10.1090/btran/151} {\path{doi:10.1090/btran/151}}.

\bibitem[Kno06]{Kn-semisimple-tensor-gltfq}
F.~Knop.
\newblock A construction of semisimple tensor categories.
\newblock {\em C. R. Math. Acad. Sci. Paris}, 343(1):15--18, 2006.
\newblock URL: \url{https://arxiv.org/abs/math/0605126}, \href
  {https://doi.org/10.1016/j.crma.2006.05.009}
  {\path{doi:10.1016/j.crma.2006.05.009}}.

\bibitem[Kno07]{Knop_tensor_envelopes}
F.~Knop.
\newblock Tensor envelopes of regular categories.
\newblock {\em Adv. Math.}, 214(2):571--617, 2007.
\newblock URL: \url{https://arxiv.org/abs/math/0610552}, \href
  {https://doi.org/10.1016/j.aim.2007.03.001}
  {\path{doi:10.1016/j.aim.2007.03.001}}.

\bibitem[Koc04]{Ko-tqfts}
J.~Kock.
\newblock {\em Frobenius algebras and 2{D} topological quantum field theories},
  volume~59 of {\em London Mathematical Society Student Texts}.
\newblock Cambridge University Press, Cambridge, 2004.

\bibitem[Kot22]{Ko}
V.~Kot\v{e}\v{s}ovec.
\newblock Asymptotics for a certain group of exponential generating functions.
\newblock 2022.
\newblock URL: \url{https://arxiv.org/abs/2207.10568}.

\bibitem[LTV23]{LaTuVa-growth-pfdim}
A.~Lacabanne, D.~Tubbenhauer, and P.~Vaz.
\newblock Asymptotics in finite monoidal categories.
\newblock {\em Proc. Amer. Math. Soc. Ser. B}, 10:398--412, 2023.
\newblock URL: \url{https://arxiv.org/abs/2307.03044}, \href
  {https://doi.org/10.1090/bproc/198} {\path{doi:10.1090/bproc/198}}.

\bibitem[LTV24]{LaTuVa-growth-pfdim-inf}
A.~Lacabanne, D.~Tubbenhauer, and P.~Vaz.
\newblock Asymptotics in infinite monoidal categories.
\newblock 2024.
\newblock URL: \url{https://arxiv.org/abs/2404.09513}.

\bibitem[LPRS24]{LaPoReSo-growth-qgroup}
A.~Lachowska, O.~Postnova, N.~Reshetikhin, and D.~Solovyev.
\newblock Tensor powers of vector representation of {$U_q(sl_2)$} at even roots
  of unity.
\newblock 2024.
\newblock URL: \url{https://arxiv.org/abs/2404.03933}.

\bibitem[Lar24]{La-char2-story}
M.~Larsen.
\newblock Bounds for {SL2}-indecomposables in tensor powers of the natural
  representation in characteristic 2.
\newblock 2024.
\newblock URL: \url{https://arxiv.org/abs/2405.16015}.

\bibitem[Mac69]{MacWilliams_orthogonal_finite_field}
J.~MacWilliams.
\newblock Orthogonal matrices over finite fields.
\newblock {\em Am. Math. Mon.}, 76:152--164, 1969.
\newblock \href {https://doi.org/10.2307/2317262} {\path{doi:10.2307/2317262}}.

\bibitem[Mat99]{Ma-hecke-schur}
A.~Mathas.
\newblock {\em Iwahori--{H}ecke algebras and {S}chur algebras of the symmetric
  group}, volume~15 of {\em University Lecture Series}.
\newblock American Mathematical Society, Providence, RI, 1999.
\newblock \href {https://doi.org/10.1090/ulect/015}
  {\path{doi:10.1090/ulect/015}}.

\bibitem[{OEI}23]{Oeis}
{OEIS Foundation Inc.}
\newblock The {O}n-{L}ine {E}ncyclopedia of {I}nteger {S}equences, 2023.
\newblock Published electronically at \url{http://oeis.org}.

\bibitem[Qua07]{calc}
J.~Quaintance.
\newblock Letter representations of {$m\times n\times p$} proper arrays.
\newblock {\em Australas. J. Combin.}, 38:289--308, 2007.
\newblock URL: \url{http://arxiv.org/abs/math/0412244v2}.

\bibitem[Soe99]{MR1698856}
W.~Soergel.
\newblock Character formulas for tilting modules over quantum groups at roots
  of one.
\newblock In {\em Current developments in mathematics, 1997 ({C}ambridge,
  {MA})}, pages 161--172. Int. Press, Boston, MA, 1999.

\bibitem[Tub24]{Tu-sandwich-cellular}
D.~Tubbenhauer.
\newblock Sandwich cellularity and a version of cell theory.
\newblock {\em Rocky Mountain J. Math.}, 54(6):1733--1773, 2024.
\newblock URL: \url{https://arxiv.org/abs/2206.06678}, \href
  {https://doi.org/10.1216/rmj.2024.54.1733}
  {\path{doi:10.1216/rmj.2024.54.1733}}.

\bibitem[TV23]{TuVa-handlebody}
D.~Tubbenhauer and P.~Vaz.
\newblock Handlebody diagram algebras.
\newblock {\em Rev. Mat. Iberoam.}, 39(3):845--896, 2023.
\newblock URL: \url{https://arxiv.org/abs/2105.07049}, \href
  {https://doi.org/10.4171/rmi/1356} {\path{doi:10.4171/rmi/1356}}.

\end{thebibliography}
\end{document}